\documentclass[14pt]{extarticle}
\usepackage[margin=1in]{geometry}
\usepackage{dsfont}
\usepackage{amsfonts}
\usepackage{tikz}
\usetikzlibrary{matrix}

\newtheorem{theorem}{Theorem}[section]
\newtheorem{pro}{Proposition}[section]
\newtheorem{lemma}{Lemma}[section]
\newtheorem{definition}{Definition}[section]

\newcommand{\proof}[1]{\noindent{\it\bf Proof:#1\ }}
\newcommand{\QED}{\hfill$\Box$\medskip}

\begin{document}

\title{   Higher-degree Smoothness of Perturbations  III}
\author{
Gang Liu \\Department of Mathematics\\UCLA }
\date{September,  2016}
\maketitle

 \section{  Introduction}
 This paper is the  sequel to \cite{5, 6}. 
 The main  results  of this paper   are the following theorems (see the relevant definitions in the later sections).
 \begin{theorem}
 	Assume that $k>1$ and $p$ is an even integer greater than $2$. Choose an positive $\nu$ with   $\nu<(p-2)/p <1$. Then there is a topological embedding 
 	${ {\cal B} }^{\nu(a)}\rightarrow { {\cal B} }^{\nu}_{C}.$
 \end{theorem}

  \begin{theorem}
  	Assume that $p$ is a positive even integer. The $p$-th power of $L_{k, \nu(a)}^p$ and $L_{k;\nu}^p$ norms, denoted by $N_{k, p, \nu(a)}:
  	W^{\nu(a)}(f, {\bf H}_{f_0})\rightarrow {\bf R} $ and $N_{k, p, \nu}:
  	V^{\nu}_{C(c^0)}(f, {\bf H}_{f_0})\rightarrow {\bf R} $  are of stratified  $C^{m_0}$ smooth in the sense that they are stratified  $C^{m_0}$ smooth 
  	viewed 	in  any slices. Moreover for $N_{k, p, \nu(a)}:
  	W^{\nu(a)}(f, {\bf H}_{f_0})\rightarrow {\bf R} $, it is stratified  $C^{m_0}$ smooth with  respect to  any of the two stratified smooth structures. Furthermore, under the assumption  that  $p>2$ and   $0<\nu<(p-2)/p$, let $I: W^{\nu(a)}(f, {\bf H}_{f_0})\rightarrow  V^{\nu}_{C(c^0)}(f, {\bf H}_{f_0})$ be the topological embedding. Then 
  	$I^*(N_{k, p, \nu}):W^{\nu(a)}(f, {\bf H}_{f_0})\rightarrow  {\bf R} $ is 
  	stratified  $C^{m_0}$ smooth in the same sense above.
  \end{theorem}
  
  Recall that $m_0=[k-2/p]$. 
  
  The immediate consequence of the above theorems is that for any  given $W_{\epsilon}^{\nu(a)}(f_0, {\bf H}_{f_0})$ considered as a local slice
   of the space ${\cal B}^{\nu(a)}$ of unparametrized stable maps,
     there 
   is a stratified $C^{m_0}$ cut-off function $\gamma_{f_0}$ such that  (i)  it is  equal to  $1$ on a slightly smaller local slice $W_{\epsilon'}^{\nu(a)}(f_0, {\bf H}_{f_0})$  with  $\epsilon'<\epsilon$,  and (ii)  it is globally defined considered as a "function" on ${\cal B}^{\nu(a)}$ although it is not supported inside $W_{\epsilon}^{\nu(a)}(f_0, {\bf H}_{f_0})$.   Indeed, let $W_{\epsilon}^{\nu(a)}(f_0, {\bf H}_{f_0})\subset V^{\nu}_{C(c^0), \delta}(f_0, {\bf H}_{f_0})$  be an embedding in the first theorem.  When  $\epsilon$ is sufficiently small, the proof of the first theorem implies that there is a
   positive constant $C$ such that for $\delta>\epsilon/C$, the desired embedding already exists. Moreover for $\delta'<\delta$ but sufficiently close to $\delta$, the above embedding induces the corresponding embedding  of the smaller neighborhoods given by  replacing $\epsilon$ by $\epsilon'$ and $\delta$ by $\delta'$.  Now let ${\hat \gamma}_{f_0}$ be a stratified $C^{m_0}$-smooth cut-off function that is  supported on  $V^{\nu}_{C(c^0), \delta}(f_0, {\bf H}_{f_0})$, hence global on ${\cal B}^{\nu}_{C}$, and is equal to 1 on $V^{\nu}_{C(c^0), \delta'}(f_0, {\bf H}_{f_0})$. Using the embedding,  the   desired cut-off function $\gamma_{f_0}$ then is the restriction of ${\hat \gamma}_{f_0}$ to $W_{\epsilon}^{\nu(a)}(f_0, {\bf H}_{f_0}).$

   Now recall the  following main result in \cite{5}.
   
   \begin{theorem}
   Let ${ W}^{\nu(a)}(f, {\bf H}_f)$ be a local uniformizer
   	(slice) centerd at $f=f_0:\Sigma_0\rightarrow M$ of the space ${\cal B}^{\nu(a)}$ of stable $L^p_{k, \nu(a)}$-maps with domains given by  ${\cal S}_t, t\in N(\Sigma_0)$. Here $N(\Sigma_0)$ is a local chart of ${\overline {\cal M}}_{0, k}$ with 
   	$t=0$ corresponding to  $\Sigma_0={\cal S}|_{t=0}$ of the central fiber of  the local universal family of stable curves ${\cal S}\rightarrow N(\Sigma_0)$.	Let $K\simeq K_t$ with $t\in N(\Sigma_0)$ be the fixed part of the local universal family, 
    and ${\widetilde W}(f_{K})$ be the corresponding space of $L_k^p$ maps with domain $K$ with associated bundle ${\cal L}^{K}\rightarrow {\widetilde W}(f_{K})$. 
   	 Then any section ${\xi}^{K}:{\widetilde W}(f_{K})\rightarrow {\cal L}^{K}$  satisfying the condition $C_1$ and $C_2$  gives rise   a stratified  $C^{m_0}$-smooth section   $\xi: { W}^{\nu(a)}(f, {\bf H}_f)\rightarrow {\cal L}$, which is still  stratified $C^{m_0}$-smooth  viewed  in any other  local slices in the sense that it is so   on their "common intersections" (=the fiber product over the space of unparametrized stable maps). 
   	
   \end{theorem}
   
    The above results together imply the existence of the desired global perturbations.
   
   \begin{theorem}
   	  For the   $C^{m_0}$-smooth section   $\xi: { W}_{\epsilon}^{\nu(a)}(f_0, {\bf H}_{f_0})\rightarrow {\cal L}$   obtained from ${\xi}^{K}$ in above theorem, the section $\gamma_{f_0}\xi$  can  be considered as    a stratified $C^{m_0}$-smooth global section on ${\cal B}^{\nu(a)} $ in the sense that  is still  stratified $C^{m_0}$-smooth  viewed  in any other  local slices. Moreover $\gamma_{f_0}\xi=\xi$  on a smaller but  prescribed 
   	   local slice ${ W}_{\epsilon'}^{\nu(a)}(f_0, {\bf H}_{f_0})$.
   	
   \end{theorem}
   
   \proof

     Clearly the section ${\xi}^{K}$ also gives rise a corresponding section, still denote by ${\xi}$,  on the larger neighborhood $V^{\nu}_{C(c^0), \delta}(f_0, {\bf H}_{f_0})$ so that $\gamma_{f_0}\xi$ is globally defined on ${\cal B}_C^{\nu}$. Hence it is globally defined on ${\cal B}^{\nu(a)}$ as well by embedding theorem. The rest of the proof is  clear.

     \QED

 	
 An immediate corollary is the existence of the   perturbed moduli spaces of stable maps as a topological manifolds with expected dimensions under the assumption that there is no non-trivial isotropy group.
 
 \begin{theorem}
 	Let $s:{\cal B}_{k, p}(A)\rightarrow {\cal L}_{k-1, p}$ be the  ${\bar{\partial }}_J $-section  of the "bundle"  ${\cal L}_{k-1, p}$ over the space of genus zero unparametrized 
 	stable $L_k^p$-maps of class $A\in H_2(M, {\bf Z})$, Here $(M, \omega)$ is    a  compact symplectic manifold 
 	with an ${\omega}$-compatible almost complex structure $J$, and 
 	the fiber of ${\cal L}_{k-1, p}$ at $(f: \Sigma \rightarrow M)$ is the space $L_{k-1}^p(\Sigma , \Lambda^{0, 1}(f^*(TM))).$ Let ${\cal W}=\cup_{i\in I} W_i$ be a sufficiently fine finite covering of the compactified  moduli space ${\overline {\cal M}}(J, A)=s^{-1}(0)$ of $J$-holomorphic sphere of class $A$,  where  $W_i=:W_{f_i}, [f]_i\in {\overline {\cal M}}(J, A)$  is a local uniformizer of a neighborhood of $[f_i]\in{\cal B}_{k, p}(A)$. Assume further  that
 	all isotropy groups are trivial and that the virtual dimension of  ${\cal M}(J, A)$ is less than $m_0$. Then there are generic small perturbations ${\nu}=\{\nu_i, i\in I\}$, which are compatible stratified $C^{m_0}$-smooth sections
 	of the local bundles $ {\cal L}_i\rightarrow W_i$, such that the perturbed moduli space ${\cal M}^{\nu}(J, A)=\cup_{i\in I}(s+\nu_i)^{-1}(0)$ is a compact  topological manifold  with  expected dimension.
 	
 \end{theorem}
 
 More generally,  combing  the theorem 1.4  above with the construction in Sec.4 of \cite {7} , we get the following  theorem that establishes the existence of the genus zero virtual moduli cycles in GW-theory.
 
 \begin{theorem}
 	Under the same assumption but without  the restriction on the isotropy groups, there are generic small perturbations ${\nu}=\{\nu_i, i\in I\}$, which are compatible sections
 	of the local bundles $ {\cal L}_i\rightarrow W_i$, such that the resulting    functorial system of  weighted perturbed moduli spaces,  $\{ (1/\#({\Gamma_I}))\cdot {\cal M}^{I, \nu}(J, A), I\in  {\cal I}\}$ is a weighted compact  topological orbifold  with  expected dimension.
 	
 \end{theorem}
  
   We refer the readers to \cite{7} for the notations used in above theorem, where the weighted  orbifold here  is called the virtual moduli cycle.

   By generalizing  the weakly smooth structure in \cite{10} to the case by allowing changes
    of the topological types of the domains as in \cite{6} and this paper,  we will prove in \cite{11} that the weighted topological orbifold ${\cal M}^{\nu}(J, A)$ is in fact 
    stratified smooth of class $C^{m_0}.$


 The first two theorems  are  proved in Sec.2 and Sec.3. 
 

 
 

 	
 For the basic facts on Sobolev spaces and calculus of Banach spaces, we refer to
  \cite{4,9}.

  \section{Topology of the space ${\cal B}^{\nu(a)}$  and ${\cal B}^{\nu}_C$}

\subsection{  The space ${\cal B}^{\nu}_C$}
  
  Let   $\Sigma=(S; {\bf p})$   be  a genus zero stable curve  modeled on a tree $T=T_{\Sigma}$  with the underlying nodal curve $S$ and (the set of) distinguished points ${\bf p}$. 
  Assume that each component $S_v, v\in T$ is marked in the sense that  a biholomorphic identification $S_v\simeq S^2$ is already given for each $v\in T$.
Then using the marking,   at each double point $d_{vu}\in S_v$  fix a small disk $D(d_{vu})=: D_{\epsilon_2}(d_{vu})$ of radius $\epsilon_2$ and an identification 
 $D^*(d_{vu})=: D(d_{vu})\setminus \{d_{vu}\}\simeq [0, \infty)\times S^1$ so that $D^*(d_{vu})$ becomes a cylindrical end with coordinate $(t, s)\in [0, \infty)\times S^1.$
  
    Now  assume that $f:S\setminus {\bf d}\rightarrow M$  is  a $L_{k, loc}^p$-map with $k-2/p>2$ that can be  continuous extended over $S$.  Given a section $\xi:S\setminus {\bf d}\rightarrow f^*TM$, its exponentially weighted $L_k^p$-norm with weight ${\nu}$ with $0<\nu<1$ is defined to be 
  $\|\xi\|_{k, p, \nu}=
  \|\xi|_{S\setminus \{ \cup_{d_{uv}} D(d_{vu})\}}\|_{k, p}+
  \Sigma_{d_{uv}\in {\bf d}}\|\xi|_{D^*(d_{vu})}\|_{k, p, \nu}.$ 
  	Here $\|\xi|_{D^*(d_{vu})}\|_{k, p, \nu}$ is defined to be $\| e_{ \nu}\cdot {\hat \xi}_{uv}\|_{k ,p}$, where ${\hat \xi}_{uv}$ is the section  corresponding to  the restriction $\xi|_{D^*(d_{vu})}$ under  the identification 	$D^*(d_{vu})\simeq [0, \infty)\times S^1$. Here 
   $e_{ \nu}(s, t)$ is a function on the   cylindrical ends defined to be:  $e_{ \nu}(s, t)=exp\{ \nu|s|\}$ on the part of cylindrical ends  with  $(s,t)\in [1, +\infty)
  \times S^1$, and $e_{ \nu}(s, t)=1 $ on $[0, 1/2]
  \times S^1$.

  Now assume that  two such continuous  $L_{k, loc}^p$-maps $f_1:\Sigma_1\rightarrow M$ and  $f_2:\Sigma_2\rightarrow M$ are in the same equivalence class so that there is
  a biholomorphic identification $\phi:\Sigma_1\rightarrow \Sigma_2$ such that 
  $f_1=f_2\circ \phi.$ In general $\phi$ does not preserve the small disks  (with  natural metrics) at double  points described above. However,  a direct computation show that for any fixed automorphism $\phi:\Sigma\rightarrow \Sigma$ and any $\xi$ above, there exist a constant $C_1=C_1(\phi) $ and  $C_2=C_2(\phi) $  such that $\|\phi^*(\xi)\|_{k, p, \nu}\leq C_1\cdot  \|\xi\|_{k, p, \nu}$ and $\|\xi\|_{k, p, \nu}\leq C_2\cdot \|\phi^*(\xi)\|_{k, p, \nu}.$
  
  Hence the notion of  exponentially decay $L_k^p$-sections $\xi$ of the bundle $f^*TM$ along cylindrical ends with weight $\nu$ is well-defined  so that  the condition $\|\xi\|_{k, p, \nu}<\infty$ is invariant with respect to the choices of the identifications of the  punctured small disks  with $[0, \infty)\times S^1$ and  the representatives $f'$ in the equivalence class $[f]$. This justifies
   the following definition.
   
   \begin{definition}
   	Let $n(d)=n(d, S)$ be the number of the (pairs of ) double points of a node curve $S$.    A continuous node  map $f:S\rightarrow M$ is said to be a $L_{k, \nu}^p$-map with asymptotic limit $c\in M^{n({\bf d})}$ if $\|"f-c"\|_{k, p, \nu}<\infty.$
   \end{definition}
   
     Here  $"f-c"$ is defined by: on the fixed part $K$, $"f-c"=f$,  and  on $D(d_{vu})$, it is equal to $f-c_{vu}$.
    Note that here  like $f$,  $c_{vu}$  is  considered as a map from $D(d_{vu})$ or $K$ to $M\subset {\bf R}^M$  so that $f-c_{vu}$ is well-defined.

  Thus such a $L_{k, \nu}^p$-map approaches its asymptotic limit $c$ along the  cylindrical ends at double points
  exponentially  with weight $\nu.$ The discussion before implies that this notion is independent of the choices of the small disks  and  the associated end structures, and that the equivalence  classes of such maps are well-defined.
  
  Now for a fixed $T$ with $n(T)=n({\bf d})$ and $c=:c_T\in M^{n(T)},$
   the space of the stable  $L_{k, \nu}^p$-maps of type $T$ with a fixed asymptotic limit 
  $c$ is defined to be 

$${\widetilde {\cal B}}^{\nu, {T}}_{c}=:
 {\widetilde {\cal B}}^{\nu, {T}}_{c_T}=\{g\,\, |\, g:S\rightarrow M,\, T(S)=T, \, \|g-c\|_{k, p;\nu}<\infty \}.$$ Here $T(S)$ is the tree associated with the nodal surface $S$.
  
  Let ${\widetilde {\cal B}}^{\nu, {T}}_{C}=\cup_{c\in M^{n(T)}} {\widetilde {\cal B}}^{\nu, {T}}_{c}$ be the space of the  stable $L_{k, \nu}^p$-maps of type $T$ and 
  ${\widetilde {\cal B}}^{\nu}_{C}=\cup_{T} {\widetilde {\cal B}}^{\nu, {T}}_{C}$
  be the space of the stable  $L_{k, \nu}^p$-maps  of all types.
  
  The  corresponding spaces of the  unparametrized  stable $L_{k, \nu}^p$-maps  will be denoted by ${ {\cal B}}^{\nu, {T}}_{c}$, ${ {\cal B}}^{\nu, {T}}_{C}$ and 
  ${ {\cal B}}^{\nu}_{C}$ respectively.
  
  Like the ordinary stable $L_k^p$-maps, the topology and local (stratified ) smooth 
  structure on these spaces of unparametrized  stable $L_{k, \nu}^p$-maps can be defined  by using local slices in a process similar to the one  in \cite{6}.
  We now give some relevant definitions briefly assuming that readers are already familiar with the notations  in  \cite{6}.

  We start with the "base" deformation $\{f_t, t\in {\bar W}(\Sigma_0)\}$ of the initial $L_{k, \nu}^p$-map $f=f_0:S_0\rightarrow M$ with asymptotic limit $c^0\in {M}^{n({\bf d})}, $ where ${\bar W}(\Sigma_0)$ is a (family of) coordinate chart(s) of $ {\bar N}(\Sigma_0)$.  Given $f_0$, the definition of $f_t$ used in this paper is the same as the one in  \cite{6}, which we refer the readers to for the detailed definitions and notations.  
  Briefly, each $t=(b, a)$, for $a=0$, the parameter $b$ describe the local moduli of the stable curves $\Sigma_b$ that has the same type $T_0$ as the initial curve $\Sigma_0$.
  Thus it is the collection of coordinates
   of the distinguished  points ${\bf p}_b$ of   $\Sigma_b$ with respect to 
    ${\bf p}_{b=0}$ of the initial curve $\Sigma_0$ upto the obvious componentwise 
    ${\bf SL}(2, {\bf C})$ action. Then $f_b$ is defined to be $f_b=f_0\circ \lambda_b^0.$ Here  $\lambda_b^0:S_b\rightarrow S_0$ is a family of componentwise  smooth maps that send the small disks on $S_b$ centered at the double points  ${\bf d}_b$ of $S_b$ to the corresponding fixed ones on $S_0$. Since $S_b$ and $S_0$ have the same desingularization, using the "marking", we only need to construct  $\lambda_b^0$  as a family self diffeomorphisms of $S^2$.  Clearly, we can 
    require that they are holomorphic aways from the small annuli abound  the "double points" and   localized on the small disks centered at these "double points".
    Indeed they  can be defined as "localized translations" near double points in the obvious way.  The parameter $a=\{a_{uv}, (u, v)\in E(T_0)\}$ is the collection of gluing parameters and  each $a_{uv}$ is a complex parameter associated to a double point  $d_{uv}$  lying on the two components $S_{0, u} $ and $S_{0, v}$ of the initial surface
    $S_0$. The gluing of the two components $f_{-}:S_{-}\rightarrow M$   and $f_{+}:S_{+}\rightarrow M$ with domains $S_-=:S_{0, u}$ and $S_+=:S_{0, v}$ at the double point $d_{uv}=d_{vu}:=d_-=d_+$ with  the complex  parameter $a\not = 0$, denoted by ${f_-}\chi_a f_+:{S_-}\chi_a S_+\rightarrow M$,  is given by (i) on $s_{\pm}<-log |a|$, it is equal to $f_{\pm};$ and on $s_{\pm}\in [-log |a|, -log |a|+1]$, ${f_-}\chi_a f_+(s_{-}, \theta)= \beta_{[-log a, -log a+1]}(s_-)f_-(s, \theta)+(1-\beta_{[-log |a|, -log |a|+1]})(s_-)c^0,$ and 
    ${f_-}\chi_a f_+(s_{+}, \theta)$ is defined by a similarly formula  but replacing  $ \theta$ by   $ \theta+\arg a.$
     Here ${S_-}\chi_a S_+$  is defined by cut of the ends $s_{\pm}> -log |a|+1$ from $S_{\pm}$ and glue back the rest along the boundary circles with an angle twisting  $\arg a.$ The general case is done inductively so that $\chi_a (f_b)$ is defined. We define $f_t=f_{b, a}=\chi_a (f_b): S_t\rightarrow M.$
 Now the metric on $S_0$  with "flat" ends at each double point induces a  corresponding metric on $S_b$, which in turn gives a metric on $S_t$ through the gluing.  It still has a flat metric on each end. Moreover  it also has a flat metric over each neck area $N(d^s)$ $\simeq $ a finite cylinder of length $2(-log |a|+1)$ which is obtained from smoothing out  a double point $d^s$ in $S_0.$ We will call each  such  neck area $N(d^s)$  an asymptotic
 end around the central circle $S^1(d)$. Thus with  these well-defined end structure on $S_t$, 
   the weight function $\nu^t$ is well-defined along both  types of  ends.
   
  Now for a fixed $t\in {\bar W}(\Sigma_0)$, a $\epsilon$-neighborhood  of $f_t$ in ${\widetilde {\cal B}}^{\nu, t}_c$ is defined to be 
  ${\widetilde W}_{c(t),\epsilon }^{\nu, t}(f_t)=\{h_t: S_t\rightarrow M\, | \|h_t-f_t\|_{k, p, \nu^t}<\epsilon\}.$ Here we use $c(t)$ to denote the ends of $f_t$. Hence for $t\in T$ in the same stratum $c(t)=c^T$ is a constant.
  
  Note that the above condition implies that $h_t$ has the same asymptotic limits $c(t)$ as $f_t$.
  
  The  
    "full" neighborhoods  ${\widetilde W}^{\nu}_{\epsilon, c}(f_0)=\cup_{t\in {\bar W}(\Sigma_0)}{\widetilde W}^{\nu, t}_{\epsilon, c(t)}(f_t).$

  A special feature of using the  exponentially  weighted $L_{k, \nu}^p$-norms is that the initial map $f_0$ can be approximated 
  by the one that is asymptotically  constants along cylindrical ends of $f_0$ at double points.
    Indeed  given any $\epsilon>0, $
   the condition $\|f-c^0\|_{k, p, \nu}<\infty$ implies that there is a $s_0$ such 
    that $\|(f-c^0)|_{[s_0, \infty)\times S^1}\|_{k, p, \nu}<\epsilon_1<<\epsilon .$ Hence we may approximate $f$ in $L_{k, \nu}^p$-norm by a map $ {\tilde f}$ such that
${\tilde f}=f$ away from the ends with $s>s_0$ and ${\tilde f}_{vu}=c^0_{vu}$ along the end at $d_{vu}$ with $s>s_0+1.$ In fact we can simply  define, assuming that there is only one end for simplicity of notations,  $ {\tilde f}=\beta_{[s_0, s_0+1]}f+(1-\beta_{[s_0, s_0+1]})c^0$, where $ \beta_{[s_0, s_0+1]}$ is a smooth cut-off function such that $\beta_{[s_0, s_0+1]}$ is equal to  $0$ for $s>s_0+1$ and $1$ for $s<s_0.$ 

Then $$\|f-{\tilde f}\|_{k, p, \nu}\leq \|(f-{\tilde f})_{S_t\setminus [s_0,\infty)\times S^1}\|_{k, p, \nu}+\|(f-{\tilde f})|_{[s_0,\infty)\times S^1}\|_{k, p, \nu}$$ $$ = \|(1-\beta)(f-c_0)\|_{k, p, \nu}\leq \epsilon_1 \|\beta\|_{C^k}<< \epsilon/8.$$

 Assume further that the gluing parameter $|a|$ is sufficiently small such that $-log|a|>>s_0$ above,  for such $a\not= 0, $ consider  the two gluings   ${f_-}\chi_a f_+:{S_-}\chi_a S_+\rightarrow M$ and ${{\tilde f}_-}\chi_a {\tilde f}_+:{S_-}\chi_a S_+\rightarrow M$. Away from the neck area, on $S_t\setminus [-log|a|,-log|a|+1)\times S^1 $, they are equal to $f_{\pm}$ and ${\tilde f}_{\pm}$ respectively so that the $(k, p, \nu)$-norm of their difference restricted  to there is less  $\epsilon_1$. On the neck area above the condition on $|a|$ implies that ${{\tilde f}_-}\chi_a {\tilde f}_+= c_0. $
  Then a   computation similar to above implies that the same conclusion for 
the $(k, p, \nu)$-norm of their difference restricted  to this area. Hence 
$\|{{\tilde f}_-}\chi_a {\tilde f}_+-{{ f}_-}\chi_a { f}_+\|_{k, p, \nu}\leq C \epsilon_1.$

Next consider a $g=(g_{-}, g_{+})\in {\widetilde W}^{\nu, t=0}_{\epsilon, c(t=0)}(f_0)$ with $f_0=f$ above. Let $\epsilon_2=\|f-g\|_{k, p, \nu}$.  Choose $\epsilon_1$ with $\epsilon_1 << \epsilon_2<\epsilon.$   Applying above argument to $g$, we get the corresponding ${\tilde g}$ and ${{\tilde g}_-}\chi_a {\tilde g}_+.$  Then for $\epsilon_1$ sufficiently small  (a)  since $\|{\tilde g}-{ g}\||_{k, p, \nu}\leq C\dot \epsilon_1$, we may  replace $g$ by ${\tilde g}$;
(b) for $|a|$ sufficiently small specified before, on the neck area
 of $S_t$, both ${{\tilde g}_-}\chi_a {\tilde g}_+$ and ${{\tilde f}_-}\chi_a {\tilde f}_+$ are equal to $c^0$  and away from the neck area they are equal to 
 ${\tilde g}$ and ${\tilde f}$ respectively so that 
    $\|{{\tilde g}_-}\chi_a {\tilde g}_+-{{\tilde f}_-}\chi_a {\tilde f}_+\||_{k, p, \nu}$ is  less  than or equal to  $\|{\tilde g}-{\tilde f}\|_{k, p, \nu} $ ($\leq \epsilon- \epsilon_2/8$,  we may assume ). This proves that there is a small "full" $\epsilon_1$-neighborhood ${\widetilde W}^{\nu}_{\epsilon_1, c}({\tilde g})$ centered at ${\tilde g}$ (or $g$) that is  contained in the given 
 $\epsilon$-neighborhood ${\widetilde W}^{\nu}_{\epsilon, c}(f_0)$ centered at $f=f_0$.

Using the argument above inductively from lower strata to higher ones,  the following 
proposition is proved.
\begin{pro}
		The topology defined by only using asymptotic constant  maps and the corresponding 
	asymptotic constant deformations is equivalent to the usual one.
	In particular, the  usual "full" neighborhoods  ${\widetilde W}_{\epsilon, c}^{\nu}(f_0)=\cup_{t\in {\bar W}(\Sigma_0)}{\widetilde W}^{\nu, t}_{\epsilon, c(t)}(f_t)$
	generates a topology on ${\widetilde {\cal B}}_c^{\nu}$.
\end{pro}

 Thus  we have defined a topology on the total space ${\widetilde {\cal B}}^{\nu}_c$  with fixed asymptotic limits $c$. Using the corresponding  neighborhoods
  $${ W}_{\epsilon, c}^{\nu}(f_0, {\bf H}_f)=\cup_{t\in {\bar W}(\Sigma_0)}{ W}^{\nu, t}_{\epsilon, c(t)}(f_t, {\bf H}_f)$$ as local uniformizers of the corresponding space ${ {\cal B}}^{\nu}_c$ of unparametrized stable maps,   a topology on this space is defined as well.

  Now for each fixed stratum labeled by $T_1\geq T_0$ and  each  $g_{t_0}:S_{t_0}\rightarrow M\in { W}^{\nu, T_1}_{\epsilon, c}(f_0, {\bf H}_f),$ 
   let  ${ W}_{\epsilon', c}^{\nu, T_1}(g_{t_0}, {\bf H}_f)$   be a small $\epsilon'$-neighborhood centered at $g_{t_0}$ and inside ${ W}^{\nu, T_1}_{\epsilon, c}(f_0, {\bf H}_f)$. 
  

   
   As in  \cite{6}, each such neighborhood will be called of the {\bf first kind} in the full neighborhood considered as an end near  $f_0$. They will be used to verify the smoothness at  $g_{t_0}$ of a function/section  on ${ W}^{\nu, T_1}_{\epsilon, c}(f_0, {\bf H}_f).$
 
  The neighborhood ${\widetilde W}_{\epsilon, c}^{\nu, T_1}(g_{t_0}, )$ in a fixed stratum $T_1\geq T_0$ can be defined similarly.

 \textbf{Note on notations:}
 
 (1) In this paper, we use subscript $c$ to denote  asymptotic limit(s) with different but almost self-evident meanings. For instance, for the "full" neighborhoods, the subscript $c$ stands
  for the collection of all $c(t)$  of the  asymptotic limits of the base deformation $f_t$, which is the finite collection of $c_{T}$ with $T\geq T_0$ where $c_T$ is the {\bf fixed } asymptotic limit  for the stratum $T$. Thus in a fixed stratum $T$, $c$ also stands for $c_{T}$ if there is no confusion. Of course, the more explicit notation $c(t)$ stands for the asymptotic limit of the base deformation$f_t$ for a fixed $t$.

  (2) In ${\widetilde W}_{\epsilon, c}^{\nu}(f_0)$  or ${\widetilde W}_{\epsilon, c}^{\nu, T_1}(g_{t_0})$, the superscript ${\nu}$ is used to indicate that the $L_k^p$-norm is measured  using a weight function along  the ends at double points
and neck area as well.  In the case, only the weight function along ends or neck area  is used 
the superscript is change into ${\nu}(d)$ or  ${\nu}(a)$ accordingly.

We note that with $\nu$ as superscript, the object  is  a family of Banach manifolds even in a fixed stratum, while ${\nu}(d)$ indicate that  it is  Banach manifold for a fixed stratum. In the case we need to distinguish
 these two cases more explicitly, the letter "W" used in the {\bf first kind} neighborhoods will be replaced by "U" for the {\bf  second kind} below.

 With this understanding of notations,  as in  \cite{6}, for the  fixed stratum $T_1\geq T_0$ as above, we use ${ U}_{\epsilon, c}^{\nu(d_{t_0}), T_1}(g_{t_0}, {\bf H}_f)$ and ${\widetilde U}_{\epsilon, c}^{\nu (d_{t_0}), T_1}(g_{t_0})$ to denote the spaces corresponding  to ${ W}_{\epsilon, c}^{\nu, T_1}(g_{t_0}, {\bf H}_f)$ and ${\widetilde W}_{\epsilon, c}^{\nu , T_1}(g_{t_0})$, but with the weight function only along the ends on the {\bf fixed} "initial" surface $S_{t_0}$. In particular, they are  Banach manifolds.  As in \cite {6}, they will be called neighborhoods of the {\bf second kind}. Note that  for the neighborhoods of the {\bf second kind} defined in \cite {6},  the metrics  for domains are the spherical ones and the $L_k^p$-norms are just the usual ones without any weight.
 
 Since in a fixed stratum with $t$  near $t_0$, the weight
  function is uniformly bounded along the neck areas. This implies the following lemma.
  
  \begin{lemma}
  On a fixed stratum $T_1\geq T_0$, the two kinds of neighborhoods  generate the  same  topology.
  \end{lemma}
 
 In order to  define the smooth structure on  ${ W}_{\epsilon, c}^{\nu, T_1}(g_{t_0}, {\bf H}_f)$ and ${\widetilde W}_{\epsilon, c}^{\nu , T_1}(g_{t_0})$, which are  the  part of ends in  a fixed stratum,  we  introduce the product structure using the same construction in  \cite{6}. Let $\Lambda=\Lambda^{t_0}: W^{T_1}(\Sigma_{t_0})\times  { W}_{\epsilon, c}^{\nu(t_0),t_0}(g_{t_0}, {\bf H}_f)\rightarrow { W}_{\epsilon, c}^{\nu, T_1}(g_{t_0}, {\bf H}_f)$
 given by $(h_{t_0}, t)\rightarrow h_{t_0}\circ \lambda_{t}^{t_0}.$
 
 Essentially the same proof in  \cite{6} implies the following.
  
   \begin{lemma}
  The identification  $\Lambda: W^{T_1}(\Sigma_{t_0})\times  { W}_{\epsilon, c}^{\nu(t_0),t_0}(g_{t_0}, {\bf H}_f)\rightarrow { W}_{\epsilon, c}^{\nu, T_1}(g_{t_0}, {\bf H}_f)$ is compatible with the natural topology of each side and 
  induces
   a smooth structure on ${ W}_{\epsilon, c}^{\nu, T_1}(g_{t_0}, {\bf H}_f)$.

  \end{lemma}

Thus this product structure defines the smooth structure for the neighborhoods of  the first kind such as ${ W}_{\epsilon, c}^{\nu, T_1}(g_{t_0}, {\bf H}_f)$ and ${\widetilde W}_{\epsilon, c}^{\nu , T_1}(g_{t_0})$,  hence the stratified smooth structure for the full neighborhood ${ W}_{\epsilon, c}^{\nu}(f_{0}, {\bf H}_f)$ and ${\widetilde W}_{\epsilon, c}^{\nu , }(f_0)$.  Note that   the neighborhoods of  the  second kind already have smooth structure as Banach manifolds. Of course, all these neighborhoods with smooth strictures are only related by continuous transition functions.

Next we extend the discussion above to allow the asymptotic limits $c$ moving.

Let $D(c(t))=:D_{\epsilon'}(c(t))$  be a small poly-disc in $M^{n({\bf d}_t)}$ centered at $c(t)$ of radius $\epsilon,$ and $\{\Gamma_{c(t)}^c, c(t)\in D(c(t)) \}$ be a smooth family of  diffeomorphisms of $M^{n({\bf d}_t)}$ that is a "translation" on $D(c(t))$ sending    the center $c(t)$ to  a point $c\in D(c(t))$  and is the identity map outside a slightly  larger poly-disc.  Here $c(t)$ is the asymptotic
limits of the base deformation $f_t$. Recall that $n({\bf d}_t)$ is the number of the  double points (or ends) of $S_t$. 
Then above base family  extended into the family $\{f_{t, c}, t\in {\bar W}(\Sigma_0)\}, \, c\in D(c(t)) \} $ defined  by composing  $f_t$ with $\Gamma_{c(t)}^c$. 
Clearly 
 this new family  is the  deformation of $f_0$  in ${\widetilde {\cal B}}^{\nu}_{C}$,and for each fixed $c$, they are in ${\widetilde {\cal B}}^{\nu}_{c}$. 
 
  Now  for each fixed $t$, let ${\widetilde W}^{\nu, t}_{\epsilon, D(c(t))}(f_t)=\cup_{c\in D(c(t))}{\widetilde W}^{\nu, t}_{\epsilon, c}(f_{t, c})$
 be the neighborhood with asymptotic limits in $D(c(t))$.
  Then the full  neighborhood as the end near $f_0$ allowing moving asymptotic limits, denoted by ${\widetilde W}_{\epsilon, C}^{\nu}(f_0)$ is defined to be $ {\widetilde W}_{\epsilon, C}^{\nu}(f_0)=\cup_{t\in {\bar W}(\Sigma_0)}{\widetilde W}^{\nu, t}_{\epsilon, D(c(t))}(f_{t, c}).$

 Next consider the case  within a fixed stratum $T_1\geq T_0$, for a  map $g_{t_0}:S_{t_0}\rightarrow M$ with the same asymptotic limit $c(t_0)$ of $f_{t_0}$,  we have defined  a  small neighborhood of type $T_1$ centered at $g_{t_0}$ with the fixed asymptotic limit $c(t_0)$, denoted by 
 $ {\widetilde W}_{\epsilon', c(t_0)}^{\nu, T_1}(g_{t_0})=:\cup_{t\in {W}^{T_1}_{\epsilon'}(\Sigma_{t_0})} {\widetilde W}_{\epsilon', c(t_0)}^{\nu, t}(g_{t})$.  Here $g_t, t\in {W}^{T_1}_{\epsilon'}(\Sigma_{t_0})$ is the deformation of $g_{t_0}$ within the same stratum.

 The corresponding neighborhood with moving asymptotic limit $c\in D_{\epsilon'}(c(t_0))$ with in the stratum is denoted by 
 ${\widetilde W}^{\nu, T_1}_{\epsilon, C}(g_{t_0})=:\cup_{t\in {W}^{T_1}_{\epsilon'}(\Sigma_{t_0}), c\in D_{\epsilon'}(c(t_0))} {\widetilde W}_{\epsilon', c}^{\nu, t}(g_{t, c})$.
 Here $g_{t, c}$ is the deformation of $g_{t_0}$ within $T_1$ but with moving asymptotic limits
 $c\in D_{\epsilon'}(c(t_0))$.
 Note that for the fixed stratum $T_1$, $c(t_0)=c(t)$ for $t\in {W}^{T_1}_{\epsilon'}(\Sigma_{t_0})$ so that
  the  definition here for a fixed  stratum is consistent with the general definition of the full neighborhoods above.
 Now   the map  $h_{t_0}\rightarrow \Gamma^{c}_{c(t_0)}\circ  h_{t_0}$ with $c\in D_\epsilon' (c(t_0))$
  for   $h_{t_0}\in {\widetilde W}^{\nu, T_1}_{\epsilon, c(t_0)}(g_{t_0})$
   gives  rise an identification  map $\Gamma:D(c(t_0))\times {\widetilde W}^{\nu, T_1}_{\epsilon, c(t_0)}(g_{t_0})\rightarrow {\widetilde W}_{\epsilon, C}^{\nu, T_1}(g_{t_0}), $ which give a product structure for ${\widetilde W}_{\epsilon, C}^{\nu, T_1}(g_{t_0}) $.  Combining with the product structure for ${\widetilde W}^{\nu, T_1}_{\epsilon, c(t_0)}(g_{t_0})$ obtained by $\Lambda^{t_0}$, we get the 
   identification $$D(c(t_0))\times W^{T_1}(\Sigma_{t_0})\times {\widetilde W}^{\nu,t_0,  T_1}_{\epsilon, c(t_0)}(g_{t_0})\rightarrow {\widetilde W}_{\epsilon, C}^{\nu, T_1}(g_{t_0}) .$$ 

Since the diffeomorphisms $\Gamma_{c}^{c(t_0)}$ only acts on the target manifolds, it is easy to see that the corresponding statement of the previous lemma still holds.

\begin{lemma}
The identification $$D(c(t_0))\times W^{T_1}(\Sigma_{t_0})\times {\widetilde W}^{\nu,t_0,  T_1}_{\epsilon, c(t_0)}(g_{t_0})\rightarrow {\widetilde W}_{\epsilon, C}^{\nu, T_1}(g_{t_0}) $$  induces a smooth structure on ${\widetilde W}_{\epsilon, C}^{\nu, T_1}(g_{t_0})$,  as a  neighborhood with moving asymptotic limits. Moreover this identification is    compatible  with  the topological structures on both sides.
\end{lemma}

Now go back the general case. We want show that these new full neighborhoods  generate a topology.  To this end, note  that  if $\|h_{t_0}-g_{t_0}\|_{\nu^{t_0}, k, p}<\epsilon$,  for 
$t\in {\bar W}_{\epsilon'}(\Sigma_{t_0})$ with $\epsilon'<< \epsilon$ and any $c\in D{\epsilon'}(c(t))$
 $\|h_{t, c}-g_{t,c}\|_{\nu^{t_0}, k, p}<\epsilon$ still holds,  where  $h_{t, c}$ and $g_{t,c}$ is defined in the same way  as $f_{t, c}$ by composing the deformation such as $h_t$ with $\Gamma_{c(t)}^c, c\in D_{\epsilon'}(c(t))$. Indeed  by the  definition of $\Gamma_{c(t)}^c, c\in D_{\epsilon'}(c(t))$, it is a smooth family in $c$ which tends to identity map in $C^{\infty}$-topology  as $c\rightarrow c(t)$ for $c\in D_{\epsilon'}(c(t))$ within any fixed stratum;  it is compatible with  process of moving from a lower stratum  to a higher  one  in the sense that, the map on the poly-disk centered at a asymptotic limit
  of the  higher stratum is the same as the one on the lower stratum.

 As  before this  is the key step to prove the following lemma.
 
 \begin{lemma}
 These full neighborhoods with moving asymptotic limits generate a topology.
 \end{lemma}

The same discussion is applicable to ${ W}_{\epsilon, C}^{\nu, T_1}(g_{t_0}, {\bf H}_f)$.  Hence we have defined  topological structures as well as   weak stratified smooth structures
 (see [L?] for the definition) on  ${\cal B}^{\nu}_C$ and ${\widetilde {\cal B}}^{\nu}_C$.

We note that as moving from a lower stratum to higher ones the number of ends are changing
in the deformation $f_t$ and $f_{t, c}$. Despite of this seemly "jumping" discontinuity, these full neighborhoods
 still generate a topology. However, it is desirable to have all the constructions above 
  with the same number of ends and asymptotic limits by introducing
   new  asymptotic ends in the higher strata. This is done by consider the  "middle" circles of the neck areas as new  asymptotic ends and requiring a further condition of fixing asymptotic limits on the construction before.
The  construction below  will be useful for later sections of this paper as well as the companion papers.


  
 Recall that  for each $t\in W^{T_1}(\Sigma_0)$, each end of $S_t$ comes  from an end of $\Sigma_b$ that  has a form $D(d_{vu})\subset S_{b, v}, v\in T_0, \, (vu)\in E(T_0)$ with $D^*(d_{vu})\simeq [0,\infty)\times S^1. $  Each such double point $d_{vu}$ above remains to be a double point in $\Sigma_t$.
  Let ${\bf d}^s_{T_1}$ be the set of double points on ${\hat S}_0={\hat S}_b$ that are smoothing out in the gluing to obtain $S_t$. Then for each $d_{uv}\in {\bf d}^s_{T_1}, $
   let $N(d_{vu}, t)\simeq [-|log|a|\,|, |log|a|\,|]\times S^1$ be a neck area of length $2|log a|$ in $S_t$ with $t\in W^{T_1}(\Sigma_{_0})$ and $S^1(d_{vu}, t)$ be the corresponding  "middle circle" with respect to the above finite cylindrical coordinate of $N(d_{vu}, t).$
    Then  $S^1(d_{vu}, t)$ will be considered as an asymptotic ends so that the number of such ends is the same as the cardinality of ${\bf d}^s_{T_1}$, and the number of total ends in the stratum $T_1$ is the same as the one of the lowest stratum $T_0.$
    
  Note that by the definition of  the deformation $f_t,$ $f_t(S^1(d_{vu}, t))=f_0(d_{vu})$ which is the asymptotic limit $c^0_{uv}$  of $f_0$ at the end at  $d_{uv}\in  {\bf d}^s_{T_1}.$ Hence the total asymptotic  limits  of deformation $f_t$
 counting the ones from middle circles  in any stratum is the same as $f_0$.
  
  Then 
  $V_{c^0}^{\nu, t}(f,{\bf H}_f)=:V_{c^0}^{\nu, t}(f_t,{\bf H}_f)$ is defined to be the collection of $g_t:(S_t, {\bf x}_t)\rightarrow (M, {\bf H}_f)$ such that  (1 )$\|g_t-f_{t, c}\|_{k, p;\nu^t}<\epsilon$; (2) the value  at $g_t$ of the  evaluation map $Ev_{vu, t}$ along  the middle circle $S^1(d_{vu}, t)$ for $d_{vu}\in {\bf d}_{T_1}^s$ 
  satisfies the condition $ Ev_{vu, t}(g_t)=\int_{S^1(d_{vu})}g_t=c^0_{vu}$. Here (i) near  each  end or neck area labeled by $d_{uv}$, we consider $g_t$ as a map to a local chart of $M$ based at $c_{u, v},$ and the values of the integrations are interpreted using such charts. (ii)   $c^0=\{c^0_{uv}, (uv)\in E(T_0)\}$ are the  asymptotic limits of $f_0$ that are the same as the total asymptotic limits  $c^0(t)$ of $f_t$, but for a $t$ in higher stratum, $c^0=c^0(t)$ decomposes as an union of $c(t)$  with the asymptotic limits
 $c({\bf d}^s_t)$ from the middle circles. both depending on $t$.
  
  Now by varying  $g_t$ it is easy the see that the above evaluation maps along middle circles   are  transversal to  the asymptotic limits $c^0(t)$ so that $V_{c^0(t)}^{\nu, t}(f_t,{\bf H}_f)$ is a smooth submanifold of $W_{c(t)}^{\nu, t}(f_t,{\bf H}_f)$  with $c^0(t)=c^0.$  Moreover  there is an identification $W_{c(t)}^{\nu, t}(f_t,{\bf H}_f)\simeq V_{c^0(t)}^{\nu, t}(f_t,{\bf H}_f)\times D_{\epsilon'}(c({\bf d}^s_t ))$
  as Banach manifolds.  
  
  Next we consider the corresponding case with moving asymptotic limits.   Denote a general point in $D_{\epsilon'}(c^0)=D_{\epsilon'}(c^0(t))$ by ${\hat c}={\hat c}(t)$ with the corresponding decomposition ${\hat c}(t)=({\hat c}({\bf d}_t), {\hat c}({\bf d}^s_t)).$ Let $f_{t, {\hat c}}=:f_{t, {\hat c}(t)}=\Gamma_{c^0(t)}^{{\hat c}(t)}\circ f_t $ be the extended  base deformation with (full) moving asymptotic limits by varying ${\hat c}(t)$. 
   Then  $V_{{\hat c}(t)}^{\nu, t}(f_{t,{\hat c}(t)},{\bf H}_f)$ is defined in the obvious way as  before so that for  each fixed stratum $T_1\geq T_0$, a neighborhood near $g_{t_0}$ with full  moving asymptotic limits is defined by   $$V_{{\hat C}}^{\nu, T_1}(g_{t_0},{\bf H}_f)=:\cup_{t\in W^{T_1}_{\epsilon'}(\Sigma_{t_0}), {\hat c}(t)\in D_{\epsilon'}(c^0(t_0))} V_{{\hat c}(t)}^{\nu, t}(g_{t,{\hat c}(t)},{\bf H}_f).$$  The full neighborhood of $f_0$ with  full  moving asymptotic limits  is  given  by  $$V_{{\hat C}}^{\nu}(f_{0},{\bf H}_f)=:\cup_{t\in {\bar W}_{\epsilon'}(\Sigma_0), {\hat c}(t)\in D_{\epsilon'}(c^0(t))} V_{{\hat c}(t)}^{\nu, t}(f_{t,{\hat c}(t)},{\bf H}_f).$$ There are similar constructions
    for 
 ${\widetilde V}_{{\hat C}}^{\nu, T_1}(g_{t_0})$ and 
   ${\widetilde V}_{{\hat C}}^{\nu}(f_{0})$, etc, accordingly.

Following  proposition will be used later in this paper.

\begin{pro}
The function $N: {\widetilde W}_c^{\nu}(f)\rightarrow {\bf R}$ given by $N(h_t)= \|h_t\|_{k, p, \nu}^p$is continuous. So is 
$N: {\widetilde W}^{\nu}(f)\rightarrow {\bf R}$.
\end{pro}

\proof

This is clear in any fixed stratum.  Arguing  inductively on strata,  we only need to  consider the case that $h_{0, a}\rightarrow h_{b=0, a=0}=h_{0, 0}$ with the gluing parameter $a\not = 0$.  Clearly we may assume that the asymptotic limit $c=0$.

Then given $\epsilon >0, $ there exists a $g_{0, 0}$ and $s_0$ such  $ |N(h_{0,0})-N(g_{0, 0})|\leq N(h_{0,0}-g_{0, 0})<\epsilon/3.$  Indeed  we can choose $g_{0, 0}$ such  that $g_{0, 0}(s, \theta)=h_{0,0}(s, \theta)$ away from the end  $s>s_0$ and $g_{0, 0}(s, \theta)=0$ for $s>s_0+1.$ Then $N(h_{0,0}-g_{0, 0})=N((h_{0,0}|_{s>s_0})<\epsilon/3$ for $s_0$ sufficiently large.

Now $|N(h_{0, a})-N(h_{0,0})\leq |N(h_{0, a})-N(g_{0,a})|+|N(g_{0, a})-N(g_{0,0})|+|N(g_{0, 0})-N(h_{0,0})|.$ For $|a|$ sufficiently small, since $g_{0, 0}=0$  for $s>s_0+1$,  so is the deformation $g_{0, a}$ so that $N(g_{0, a})=N(g_{0,0})$.   Then it follows from that fact that at the lowest stratum, 
$h_{0,0}$ is in the $2/3\epsilon $-neighborhood of $g_{0, 0}$, so is $h_{0, a}$ in the "full" $2/3\epsilon $-neighborhood of $g_{0, 0}$ for $|a|$ sufficiently small.
 Since any such a neighborhood is also one of the neighborhoods entered at $h_{0, 0}$, this implies that there is a $\delta>0,$  such that for $|a|<\delta,$  
 $N(h_{0,a}-g_{0, a})<2/3\epsilon.$  Hence $|N(h_{0, a})-N(g_{0,a})|\leq N(h_{0,a}-g_{0, a})<2/3\epsilon$ for such $a$.

  The last statement allowing the asymptotic limits $c$  moving
   can be derived from the above easily. We leave it to the readers.
   
 \QED

\subsection{ The mixed norm for ${\cal B}^{\nu(a)}$}
  
  Now  recall  that the space 
  	${\tilde {\cal B}}^{\nu(a)}=\cup_{T}{\tilde {\cal B}}^{\nu(a), T}$ is the collection  of the parametrized stable $L_{k, \nu(a)}^p$-maps in \cite{6}. In the analytic setting  for  Floer homology of  \cite{7}, it was already implicitly used.
  The topology on the ${ {\cal B}}^{\nu(a)}$ of unparametrized stable maps  is  generated  by the local slices defined in 
  \cite{6},  denoted by $W^{\nu(a)}(f, {\bf H}_f)=\cup_{T\geq T_0} {W}^{\nu(a), T}(f, {\bf H}_f)$ as the union of the strata, where the center     $f=f_0:(S_0, {\bf x})\rightarrow (M,  {\bf H}_f)$  is in the lowest stratum of type $T_0$. Recall that  here the elements $g$ in  the lowest stratum ${W}^{\nu(a), T_0}(f, {\bf H}_f)$ are $L_k^p$-maps measured in the spherical metric of the domains, but for elements $g_{t}:S_{t}\rightarrow M$ with $t_0=(b, a), a\not = 0$ in the higher stratum ${W}^{\nu(a), T_1}(f, {\bf H}_f)=\cup_{t\in W^{T_1}(\Sigma_0)} {W}^{\nu(a),t,  T_1}(f, {\bf H}_f)$, an  exponentially weighted  $L_k^p$-norm is used with a wight $\nu $ (denoted by $\nu(a)$) along each neck area.
  		Thus in the definition of $W^{\nu(a), T_1}(f, {\bf H}_f)$ a "mixed" $L_k^p$/$L_{k, \nu(a)}^p$-norm is used. Unlike the $L_{k, \nu}^p$-norms 
  	before, these mixed norms appear to be "discontinuous" as  passing from higher strata to lower ones. Despite of this, they still define a topology on the space 	${ {\cal B}}^{\nu(a)}$.
    
    However,   the easiest way to prove this is to  define  a modified version of above mixed norm by rescaling. The new mixed norm  are equivalent to the old ones, 
     but  it is at least "upper continuous"/"monotonically  increasing" with respect to degenerations from high strata to the lower ones.  As  explained in next
      subsection, this monotone  property   implies that any element in a neighborhood defined using  the  scaled norm is always deformable.  Consequently  these new neighborhoods generate a topology, so do the old ones.


    

\subsection{The rescaled mixed  norm and embedding theorem}
  
 The definition of the rescaled mixed norms relies on the    the estimate and the  embedding theorem for fixed strata  in  \cite{5} that we recall  now.
   
 \begin{lemma}
 	Assume that $k>1$ and $p$ is an even integer greater than $2$. Choose an positive $\nu$ with   $\nu<(p-2)/p <1$. Then  for $\xi$
 	with $\xi(0, 0)=0$ at $(x, y)=(0, 0),$  $$||\xi||_{k, p;\nu}\leq C||\xi||_{k, p}.$$ 	
 \end{lemma}
 
  Also recall the following embedding theorem in \cite{5} for the case that the domain is the fixed $S^2$ ( stated in the notations used there).
 
 \begin{theorem}
 	
 	Under  the  assumption in the above lemma, there is a smooth embedding ${\widetilde {\cal B} }_{k,p}(c)\rightarrow {\widetilde{\cal B} }_{k,p;\nu}(c)$. Here ${\widetilde {\cal B} }_{k,p}(c)\subset {\widetilde {\cal B} }_{k,p}$ is the collection of $L_k^p$-maps
 	$f:(\Sigma, x)\rightarrow (M, c)$ from $\Sigma$ with one marked point $x$ to $M$ sending $x$ to a fixed constant $c\in M,$ and ${\widetilde {\cal B} }_{k,p;\nu}(c)$ is defined similarly with $f$  approaching to $c$ exponentially with weight $\nu$ along the cylindrical  end at $x$.

 \end{theorem}

  With  above preparations, we can defined the the rescaled mixed norms  and at the same time 
 prove  the main theorem of this section.
  
 \begin{theorem}
 Assume that $k>1$ and $p$ is an even integer greater than $2$. Choose an positive $\nu$ with   $\nu<(p-2)/p <1$.
 Then  after rescaling of the relevant norms given below,  the uniformizers $W^{\nu(a)}(f, {\bf H}_f)$, $[f]\in  {\cal B}^{\nu(a)}$  defines  a topology on ${ {\cal B}}^{\nu(a)}$. Here $f\in [f]$ is a representative
    of the unparametrized stable map $[f].$
Consequently there is a topological embedding 
 ${ {\cal B} }^{\nu(a)}\rightarrow { {\cal B} }^{\nu}_{C}.$
 \end{theorem}
  
  \proof

 To prove the first part of the theorem, we only need to show that, for any given 
  element $g_{t_0}\in W^{\nu(a), T_1}_{\epsilon }(f, {\bf H}_f)$ in the stratum $T_1$, 
  there is a deformation $g_t$ inside 
  $W^{\nu(a)}_{\epsilon }(f, {\bf H}_f)$.  Indeed assume that  this is true, then for $t$  sufficiently close to $t_0$ and $\epsilon'<<\epsilon$, by triangle inequality applying to  the $t$-dependent family of norms, 
  the collection of stable maps $h_t$ with $||h_t-g_t||_{k, p, \nu(a)}<\epsilon'$, $t\in {\bar W}(\Sigma_{t_0})$, denoted by $W^{\nu(a)}_{\epsilon' }(g_{t_0}, {\bf H}_f)$ is contained inside 
  $W^{\nu(a)}_{\epsilon }(f, {\bf H}_f)$.  Hence  these neighborhoods form a topological basis.
 
 
  To prove  the existence of the  above deformation, 
  we  assume that the mixed norm is already monotone increasing  during the degenerations first. Then for any given $g_{t_0}\in W^{\nu(a), T_1}_{\epsilon }(f, {\bf H}_f)$ in the stratum $T_1$, 
  let  $g_t$ with $|t-t_0|<\epsilon_1$ be the (any) deformation constructed before.
  Then for $\epsilon_1$ sufficiently small, by monotonicity above  $\|g_t\|_{k, p, \nu(a)}\leq \|g_{t_0}\|_{k, p, \nu(a)}<\epsilon
  $ so that the deformation $g_t$ is inside $g_{t_0}\in W^{\nu(a), }_{\epsilon }(f, {\bf H}_f)$.
     
  The desired monotonicity/up-semi-continuity follows if the embedding constant  $C$ of  the last lemma is equal to 1.

  To see this,  assume that as $t\rightarrow t_0$ from a higher stratum
   labeled by $T_1$ to a lower one of $T_0$, the stable map $h_t:S_t\rightarrow M $ is a degeneration to  $h_{t_0}: S_{t_0}\rightarrow M$  constructed before with $h_t$ and $h_{t_0}$ in $W^{\nu(a)}_{\epsilon}(f, {\bf H}_f)$.  Clearly by induction, we only need to consider the essential case that $h_t$ only has neck ares without any other ends (i.e. in the top stratum) while $h_0$ is in the lowest stratum.
   By the embedding theorem in  \cite{6} recalled above,  both of then are in $W^{\nu}_{\epsilon}(f, {\bf H}_f)$ with $h_t\rightarrow h_{t_0}$ as $t\rightarrow t_0$ in $L_{k,p, \nu}$-topology. 
  
   Now the continuity of $L_{k,p, \nu}$-norm proved in previous section implies that $\|h_{t_0}\|_{\nu, k,  p}=\lim_{t\rightarrow t_0}\|h_{t}\|_{\nu, k,  p}.$ But by assumption that $T_1$ is the top stratum with only one neck area so that  $\|h_{t}\|_{\nu, k,  p}=\|h_{t}\|_{\nu(a), k,  p}$. Hence the embedding theorem with assumption that $C=1$ in \cite{6}  implies that 
  $\|h_{t_0}\|_{\nu(a), k,  p}\geq \|h_{t_0}\|_{\nu, k,  p}=
  lim_{t\rightarrow t_0}\|h_{t}\|_{\nu(a), k,  p}.$
  
  This finishes the proof of  the theorem under the assumption that $C=1$.
 
To get rid of above assumption,   we   define below a new equivalent mixed norm 
 that is monotone increasing
during degenerations and hence    prove the theorem.

\medskip

\noindent   
{\bf $\bullet$  Monotonicity/Up-semi-continuity of the  Mixed Norm}
 
 The desired new $L_{k, \nu(a)}^p$-norm is defined by  using a  rescaling the metric of the domains.
 Let $S_t$ be the underlying nodal curve with $t=(b, a)$ so that it is obtained by gluing with the parameter $a$ from the  curve $S_b$ with $b\in W^{T_0}(\Sigma_0)$ in the lowest stratum. Then on   the fixed  part $K_t=K_b=K_0$, the  rescaled  metric is the same as the spherical one on $K_b$ and $K_0$,  and on all the neck areas, it is still the cylindrical  one  defined before. On all the small disks $D({d_{uv, t}})=:D({d_{uv}})$ centered at the double point ${d_{uv, t}}$, the new metric is the spherical (or the flat) one on  $D({d_{uv}})$ before multiplying  by the constant $C$ in the estimate in above lemma. Here $d_{uv, t}$ is  one of the double points that comes from the corresponding  double point $d_{uv}$ of  $S_0$. 
 Using  cut-off functions supported in an annulus around the boundary of $D({d_{uv, t}})$,   we get a family of desired smooth metrics parametrized by $t$ in the sense that (i)   it is a smooth family within each  stratum; (ii) as the parameters $t$ in  a higher stratum tends to $t_0$ in a lower stratum,  each of the finite cylindrical  metrics  along  a neck area of $S_t$  "tends" to the  cylindrical  metrics  on  the two ends ($\simeq$ the two punctured small disks) at a corresponding double point of $S_{t_0}$
 while the new metric on  these two small disks are the above rescaled spherical/flat ones. Thus  during the deformation, the new family of metrics like the old ones is not continuous  exactly along those relevant necks. 
 
 The new $L_{k, \nu(a)}^p $-norm, denote by $\|-||_{k, p, \nu^{new}(a)}$ is defined by using the new metrics on the domains. Clearly  it is equivalent to the old one.

Despite of the discontinuity of the new metric, at least  in the case with fixed  asymptotic limit $c$, in the neighborhood
 $W^{\nu^{new}(a)}(f, {\bf H}_f)$, the mixed norm $\|-||_{k, p, \nu^{new}(a)}$ is monotone decreasing by the construction based on the estimate in the previous lemma.
 
 This also proves the second half of the theorem on the embedding  at least for the case with fixed  asymptotic limit $c$. The general case can be easily derived from this case. We leave it to the readers.
 
 \QED
  \medskip

\section{Stratified smoothness of the $p$-th power of  $L_{k;\nu(a)}^p$ and $L_{k;\nu}^p$  norms}
 
 The main theorem of this section is the following.
 
 \begin{theorem}
 Assume that $p$ is a positive even integer. The $p$-th power of $L_{k, \nu(a)}^p$ and $L_{k;\nu}^p$ norms, denoted by $N_{k, p, \nu(a)}:
 W^{\nu(a)}(f, {\bf H}_{f_0})\rightarrow {\bf R} $ and $N_{k, p, \nu}:
 V^{\nu}_{C(c^0)}(f, {\bf H}_{f_0})\rightarrow {\bf R} $  are of stratified  $C^{m_0}$ smooth in the sense that they are stratified  $C^{m_0}$ smooth 
 viewed 	in  any slices. Moreover for $N_{k, p, \nu(a)}:
 	W^{\nu(a)}(f, {\bf H}_{f_0})\rightarrow {\bf R} $, it is stratified  $C^{m_0}$ smooth with  respect to  any of the two stratified smooth structures defined in \cite{6}. Furthermore, under the assumption  that  $p>2$ and   $0<\nu<(p-2)/p$, let $I: W^{\nu(a)}(f, {\bf H}_{f_0})\rightarrow  V^{\nu}_{C(c^0)}(f, {\bf H}_{f_0})$ be the topological embedding. Then 
 	$I^*(N_{k, p, \nu}):W^{\nu(a)}(f, {\bf H}_{f_0})\rightarrow  {\bf R} $ is 
 	stratified  $C^{m_0}$ smooth in the same sense above.
 \end{theorem}
 
 In the rest of this section, we use the notations and constructions
 in \cite {5, 6} without detailed explanations.
 
 To prove the theorem, we  make a few reductions.


 (1) Using the family of coordinate charts  $Exp_{f_{t, c}}$ with $t\in { N}^{T_1}(\Sigma_{t_0})$ and $c\in D(c^0(t_0))$ in the stratum $T_1$,
    $W^{\nu, T_1}_{C(c^0)}(f, {\bf H}_{f_0})$ can be identified with a family of  open balls in Banach spaces parametrized by $(t, c)\in N^{T_1}(\Sigma_{t_0})\times D(c^0(t_0)).$   Hence we simply assume that  it 
   is  a family of  Banach spaces. Similarly for $W^{\nu(a), T_1}(f, {\bf H}_{f_0}).$  
  (1) Note that the diffeomorphisms $\Gamma^{\hat c}_{{\hat c}^0}, {\hat c}\in D({\hat c}^0)$ only act on the 
  target space $M^{n({\bf d}_0)}$, similarly for $\Gamma^{ c}_{{ c}(t)}, {c}\in D({\hat c}(t))$. Then it is easy to see that the parameter $c$/${\hat c}$  appears in  the express $p$-th powers of these
  norm functions below  only as a smooth "factor" so that these functions depend smoothly on $c$/${\hat c}$.
  Hence we only need to consider the case  with fixed asymptotic limits.
  
  (3)  The key step that reduces the case here to the ones in  \cite{5} is to 
   introduce the enlarged corresponding three spaces: 
   (a) the space of $L_{k, \nu(a)}^p$-maps centered at $g_{t_0}$, ${\widehat  W}^{\nu(a), T_1}_c({\hat g}_{t_0} )=:\prod_{v\in T_1} {\widetilde W}^{\nu(a), T_1}_c({\hat g}_{t_0, v} )$, (b) the spaces of usual $L_k^p$-maps
  ${\widehat  U}^{T_1}_c({\hat g}_{u_0})=:\prod_{v\in T_1 }{\widetilde U}^{T_1}_c({\hat g}_{u_0, v})$, with $u_0=u^{T_1}(t_0)$ and $g_{u_0}=g_{t_0}\circ \phi_{t_0}$ as we did similarly in \cite{6},  and (c) the corresponding spaces of $L_{k, \nu}^p$-maps. Clearly we only need to prove the corresponding results for these spaces.
  Here ${\hat g}_{t_0}$ is the desingularization of $g_{t_0}.$
 Then   one can apply the corresponding results or arguments in \cite{5} componentwisely.

  Note that except the last statement, the proofs of the theorem for the two cases  are similar, and it appears that the proof for the case for $N_{k, p, \nu(a)}$ is harder. Hence we will only deal with this case.  Since in this case, there is no need to fix $c$/${\hat c}$, we drop this restriction and only consider the cases already in \cite{6}.

   First note that the restriction of $N_{k, p, \nu(a)}$ to any "fiber" with fixed $t\in T_1$, $N^{t, T_1}_{k, p, \nu(a)}:\prod_{u\in T_1} {\widetilde W}^{\nu(a),t,  T_1}({\hat g}_{t_0, u} )\rightarrow {\bf R}$ is smooth,  by   applying  componentwisely  the corresponding result in  \cite{5} for the   fixed domain $S^2$.  In particular, it  is smooth at the central slice
   $\prod_{u\in T_1} {\widetilde W}^{\nu(a), t_0, T_1}({\hat g}_{t_0, u} )$.
  Recall that the smooth structure on $\prod_{u\in T_1} {\widetilde W}^{\nu(a),  T_1}({\hat g}_{t_0, u} )$ is obtained by using the identification $\prod_{u\in T_1} {\widetilde W}^{\nu(a),  T_1}({\hat g}_{t_0, u} )\simeq \{\Pi_{u\in T_1} {\widetilde W}^{\nu(a),t_0,  T_1}({\hat g}_{t_0, u} )\}\times N^{T_1}(\Sigma_{t_0})$ induced by the maps $\lambda^{t_0}_t, t\in N^{T_1}(\Sigma_{t_0})$.  Thus the smoothness of $N_{k, p, \nu(a)}$ in this case  will follow from the smoothness the  composition of   maps  $\{\prod_{u\in T_1} {\widetilde W}^{\nu(a),t_0,  T_1}({\hat g}_{t_0, u} )\}\times N^{T_1}(\Sigma_{t_0})\rightarrow \Pi_{u\in T_1} {\widetilde W}^{\nu(a),  T_1}({\hat g}_{t_0, u}) $  and  $N_{k, p, \nu(a)}:\prod_{u\in T_1} {\widetilde W}^{\nu(a),  T_1}({\hat g}_{t_0, u}) \rightarrow {\bf R}$.
  
  Then  the composed  map can be  interpreted as  a function 
  on $(\xi, t)$, where  $t\in N^{T_1}(\Sigma_{t_0})$ and ${\xi}$ is  a $L_k^p$-section on the  fixed desingularization ${\hat S}_{t_0}$ but with a varying metric $h_t$ parametrized by $t$ that is  used to defined the $L_{k, \nu(a)}^p$-norm of $\xi.$

  More specifically, write $N_{k, p, \nu(a)}$ as the summation $N_{k, p, \nu(a)}=\Sigma_{j=0}^k N_{(j), p, \nu(a)}.$ For simplicity, we only consider the 
  first term of the summation, $N_{(0), p, \nu(a)}.$
  
  Let $m=p/2$. Then in the coordinate chart above, the  first  summand of the composed function, still denoted by  
  $$N_{(0), p, \nu(a)}(\xi, t)=\int_{{\hat S}_{t}}  <\xi\circ \lambda^{t_0}_t,  \xi\circ \lambda^{t_0}_t>^m\cdot e_{\nu(a),t}\cdot dvol_{h_t}
  $$ $$ =\int_{{\hat S}_{t_0}}  <\xi,  \xi>^m  e_{\nu(a),t}\circ (\lambda^{t_0}_t)^{-1}\cdot    Jac ^{-1} (\lambda^{t_0}_t)\cdot (det^{1/2} (h_t\circ h^{-1}_{t_0}))\cdot dvol_{h_{t_0}}.$$
  
  Here $ e_{\nu(a),t}$ is the weight function along the neck area on $S_t$.  The key point now is that the   metrics $h_t$ and diffeomorphisms ${\lambda^{t_0}_t}$ ( and  $\phi_t$ used latter) are  smooth in $t$ and only contribute   "weight functions" to the integrand smooth in $t$  so that $N_{(0), p, \nu(a)}$  is indeed smooth  by the same corresponding argument in  \cite{5}.
  Hence  $N_{k, p, \nu(a)}$ is smooth.

  The proof of the rest of the theorem is just  a adaption
   of the above argument in the other cases in the obvious manner. For completeness, we give the details.

  
  Similar argument shows  that  the  transformed $N_{k, p, \nu(a)}$ by $\Phi$, denoted by  ${ N}_{k, p, \nu(a)}^{\Phi}$ , defined on $ \prod_{v\in T_1 }{\widetilde U}^{T_1}({\hat g}_{u_0, v})$ is smooth as well. Here $\prod_{v\in T_1 }{\widetilde U}^{T_1}({\hat g}_{b_0, v})$  is either with the  usual smooth structure or  the "product smooth structure" given by  the identification $ \prod_{v\in T_1 }{\widetilde U}^{T_1}({\hat g}_{b_0, v})\simeq \prod_{v\in T_1 }{\widetilde U}^{u_0, T_1}({\hat g}_{u_0, v})\times N^{T_1}(\Sigma_{u_0})$.
  
  Indeed the map $\Phi$ is induced by the identification maps $\phi_t, t\in N^{T_1}(t_0)$  from $S_{u_0}$     to $S_{t}$. 
  After conjugating with $\phi_t$, the diffeomorphisms  $\lambda^{t_0}_t$ become $\gamma^{u_0}_u$
   with $u=u^{T_1}(t)$. Hence  with respect to  the  "product smooth structure", the first summand of ${ N}_{k, p, \nu(a)}^{\Phi}$ defined on ${\widehat  U}^{T_1}({\hat  g}_{u_0})$, evaluated at  $(\eta, u)\in  \prod_{u\in T_1 }{\widetilde U}^{u_0, T_1}({\hat g}_{u_0, v})\times N^{T_1}(\Sigma_{u_0})$, becomes
   $$N^{\Phi}_{(0), p, \nu(a)}(\eta, u)=\int_{{\hat S}_{u}}  <\eta\circ \phi^{-1}_t,  \eta\circ \phi^{-1}_t>^m e_{\nu(a),t}\cdot dvol_{h_t}
   $$
   $$ =\int_{{\hat S}_{u_0}}  <\eta,  \eta>^m e_{\nu(a),t}\circ \phi_t \cdot \phi_t^*(dvol_{h_t}).
   $$
    Since $t=t(u)$ is smooth in $u$, so are  $\phi_t$ and $ \phi_t^*(dvol_{h_t})$,  this proves that $N^{\Phi}_{(0), p, \nu(a)}(\eta, u)$ is smooth in $(\eta, u)$ for the same reason as the previous one.
    
    As for ${\widehat  U}^{T_1}({\hat g}_{u_0})$  with usual smooth structure, we note that in this case, no identifications of the  domains  by  non-holomorphic $\gamma_{u}^{u_0}$ are needed; the underlying surface of the domain $(\Sigma_u, {\bf p}_u), u\in N^{T_1}(\Sigma_{u_0})$ is the fixed ${\hat S}_{u_0}$ with the distinguished points  ${\bf p}$ moving whose positions are parametrized by $u$.
    Of course, the identification  $\phi_t:{\hat S}_{u_0}\rightarrow {\hat S}_{t}$ is the same as before. Then the first summand of ${ N}_{k, p, \nu(a)}^{\Phi}$ defined on ${\widehat  U}^{T_1}({\hat g}_{u_0})$ , evaluated at  $\xi$, formally takes the same  form, 
    
    $$N^{\Phi}_{(0), p, \nu(a)}(\xi)=\int_{{\hat S}_{t}}  <\xi\circ \phi^{-1}_t,  \xi\circ \phi^{-1}_t>^m\cdot e_{\nu(a),t}\cdot dvol_{h_t}
    $$
    $$ =\int_{{\hat S}_{u_0}}  <\xi,  \xi>^m \cdot e_{\nu(a),t}\circ \phi_t \cdot \phi_t^*( dvol_{h_t}).
$$
  
   Now in this smooth structure  of ${\widehat  U}^{T_1}({\hat g}_{u_0})$, $u$ and $t$   are  smooth function of $\xi$. So is $\phi_t$ so that $\phi_t^*( dvol_{h_t})$ is smooth in $\xi.$ The conclusion follows now again by the same argument in  \cite{5}. 
   

  Hence we conclude that $N_{k, p, \nu(a)}$ defined on $W^{\nu(a), T_1}(g_{t_0}, {\bf H}_{f_0})$  and ${ N}_{k, p, \nu(a)}^{\Phi}$ defined on ${ U}^{T_1}({ g}_{b_0}, {\bf H}_{f_0})$ with respect to the  smooth structures  above
   are of class $C^{m_0}$.  
   
   Since transition function from the given slice $W^{\nu(a)}(f, {\bf H}_{f_0})$   above to another slice is induced by reparametrizations $\psi=\psi_{\xi}, \xi\in W^{\nu(a)}(f, {\bf H}_{f_0})$ which are  $C^{m_0}$ smooth in $\xi$, a similar argument as the last one above implies that  $N^{\Phi}_{(0), p, \nu(a)}$ (or $N^{\Phi}_{(0), p, \nu})$ is  stratified smooth of class  $C^{m_0}$ viewed in any other slices.
  
  Finally to prove the last statement, we use the reduction  (2) so that we only need to consider the case  for the slice  ${ U}^{T_1}_{c}({ g}_{b_0}, {\bf H}_{f_0})$  of $L_k^p$-maps with fixed values $c$ at double points.
  
  Then the formula for first  summand of $I^*({ N}_{k, p, \nu})$ defined on ${ U}^{T_1}({ g}_{u_0}, {\bf H}_{f_0})$, evaluate at  $\xi$ is given by the same formula as the one above for $N^{\Phi}_{(0), p, \nu(a)}(\xi)$ but  replacing the weight function $e_{\nu(a),t}$ only along neck areas there  by the weight function $e_{\nu,t}$ along both  the neck areas and the cylindrical ends  at double points. More specifically 
  $$I^*({ N}_{k, p, \nu})(\xi)=\int_{{\hat S}_{t}}  <\xi\circ \phi^{-1}_t,  \xi\circ \phi^{-1}_t>^m\cdot e_{\nu,t}\cdot dvol_{h_t}
  $$
  $$ =\int_{{\hat S}_{u_0}}  <\xi,  \xi>^m \cdot e_{\nu,t}\circ \phi_t \cdot \phi_t^*( dvol_{h_t}).
$$
  
  This is a multi-linear function in $ <\xi,  \xi>^m $, $e_{\nu,t}\circ \phi_t  $ and $\phi_t^*( dvol_{h_t})$ (or $ dvol_{{\hat S}_{b_u}}$),  and these terms are smooth in $\xi$.
  	
  	The estimate in  \cite{5} stated in the last lemma of Sec. 2, implies that this function is bounded  if $\xi $ is in $L_2^p$, which is indeed the case for  $\xi $ here. Hence 
  	as in  \cite{5}, the boundedness of the   multi-linear map implies the required smoothness.

\end{document}